\documentclass[12pt]{article}
\usepackage[utf8]{inputenc}
\usepackage[T1]{fontenc}

\usepackage{palatino}

\usepackage{amsfonts}

\usepackage{mdframed}

\usepackage{amssymb, amsmath}
\usepackage[english]{babel}
\usepackage{graphicx}
\graphicspath{ {./images/} }
\usepackage{subcaption}
\usepackage[export]{adjustbox}
\usepackage{wrapfig}
\usepackage{url}
\usepackage{setspace} 

\usepackage{authblk}

\title{Andrej N. Kolmogorov's in front of the “affaire Lysenko”: an episode in mathematics and ideology in the 20th century}
\author{\textbf{Fascitiello Isabella}}
\affil{University of Roma Tre, Department of Mathematics and Physics}

\author{Fascitiello Isabella}

\date{\today}

\begin{document}
\maketitle
\tableofcontents

\begin{abstract}

Andrej N. Kolmogorov's life and outstanding mathematical work was intertwined with the political and cultural evolution of the Russian Empire from the Tzarist regime to the establishment and development of the Soviet Union. His early involvement in scientific education (and his own education) was marked by the progressive education movement continuing under the Lenin rule. Under the dramatic years of Stalin rule he was involved in the process against Nikolai N. Luzin but had a conflict in 1940 with Trochym D. Lysenko, thus founding himself in accordance and then in opposition to Ernst Kol'man (1892–1979), a very influential figure in mathematical circles. Kolmogorov's political views had been considered either that of a authentic Marxist (Graham 1993) loyal to the regime or that of a great scientist, a member of the intelligentsia science, who lived through the 1930s and 1940 acting sometimes against his principles (Arnold 2000, Lorentz 2002). 

In this paper the context of the 1940 mathematical paper by Kolmogorov on Mendel's laws of genetics is considered, also thanks to a recent study by Stark and Seneta (2011). Kolmogorov's vision of the relationship between mathematics and science in the understanding of the natural world appears to have encouraged him in tacking part in the discussions around Mendelism in while Lysenko's position against it was becoming the only permissible for Soviet Union biologists.   

\end{abstract}

\section{Kolmogorov, a representative of Russian “intelligentsia science” and his education}
The young Kolmogorov, born two years before the 1905 upheavals that marked the final years of the Tzarist regime, was raised in an environment soaked by progressive education: after some years of home schooling at his grandfather home until 1910 under the guidance of his aunts, which adhered to the principles of advanced "democratic" pedagogy then in vogue in Russia and which were the heart of school policy in the first years after the October Revolution and until 1931\footnote{See [Karp, Vogeli 2010], [Karp 2012] and [Karp, Schubring 2014]}, since the age of 7 he attended a private gymansium in Moscow, founded by two woman, Evgenja (Evgeniya) Albertovnava Repman and Vera Fedorova\footnote{See [Tikhomirov 1988]. E.A. Repman (1870-1937), founder and director of the school, was the eldest daughter of  Albert Hristianovitch Repman (1834-1917), sinced 1889 director of the section on applied physics of the Polytechnical Museum of Moscow (founded by the zar Alexander II in 1870 as Museum of Applied Knowledge). After 1917, it was renamed Section grade school no. 23. Biographical essays [Tikhomirov 1988] and [Shyraev 2000] draw many informations from Kolmogorov's reminiscences included in his posthumous book on mathematics as science and as a profession [Kolmogorov 1988] published in Russian in a series by the magazine \textit{Kvan}. In 1942 Kolmogorov would marry Anna Dmitrievna Egorova, whose first husband was the mathematician and painter Sergei Mikhailovich Ivashyov-Musatov, the mother of Oleg Sergeyevich Ivashyov-Musatov, who majored in mathematics with his stepfather. Kolmogorov, Anna Dmitrievna Egorova and her ex-husband were former students of the Repman school}. Thus, the final years of highschool and those as university student were marked by the end of Tzarist monarchy in October 1917 – he was 14 years old – and the rise of bolshevism under Lenin. He had to leave Moscow in 1918-20 with his family, as he himself recounts in [Kolmogorov 1988]: 
\\[1em]

\leftskip=1cm
\noindent
In 1918-1920 life in Moscow was not easy. In schools, only the most persistent were seriously engaged. At that moment, together with my elders, I had to leave for the construction of the Kazan-Ekaterinburg (now Sverdlovsk) railway.
\\[1em]

\leftskip=0cm
\noindent
By the cultural atmosphere of his early years, he was imbued by the  patriotic scientism of Russian radical intelligentsia, and no wonder that in 1920 he enrolled both at the Physics and Mathematics Department of Moscow University and at the Institute of Chemical Engineering “D.I. Mendeleev”. \textit{Engineering was then perceived as something more serious and necessary than pure science}, he would say in 1963 during and interview with the magazine \textit{Ogonek}. Furthermore in 1922 – he was he was enrolled as a teacher of mathematics and physics and boarding school educator in a secondary school of the network under the administration of the People's Commissariat of Education (known as Narkompros), led by Anatoly Lunacharsky (1875-1933) with Lenin's wife, Nadezhda Konstantinovna Krupskaya (1869-1939): 
\\[1em]

\leftskip=1cm
\noindent
Now I remember with great pleasure my work at the Potylikha Experimental School of the People's Commissariat of Education of the RSFSR. I taught mathematics and physics (at that time they were not afraid to entrust the teaching of two subjects to 19-year-old teachers at the same time) and took an active part in the life of the school (I was the secretary of the school board and a boarding school educator).\footnote{Although perhaps this early involvement in elementary education it began as a job out of necessity rather than will, [Abramov 2011] shows that in the 1960s and 1970s he returned to his interest in education, marked by his own experience, participating in efforts to improve the secondary mathematics education in the USSR}.
\\[1em]

\leftskip=0cm
\noindent
Alexander Karp has descrived the influence on mathematical education in the Russian Empire of the modern school movement; his description help us understanding the cultural atmosphere in which the young Kolmogorov work as mathematics and physics teacher:    
\\

\leftskip=1cm
\noindent
In place of all existing types of educational institutions, a statute of 1918 established the so-called unified labor schools. These schools were divided into two stages, and the network of first-stage schools, which were far more numerous to begin with, continued to be intensively developed. 
\\
$[...]$ The goal was to eliminate from schools anything reminiscent of former discipline and drills, including exams, textbooks, and even separate subjects (including mathematics). The ideas of American progressive educators were taken up and developed in Russia (Soviet Union); schools made use of projects, laboratory work, group work, and, above all, "complexes". 
\\
“Complexes” had to link through one overarching theme topics that had previously been studied in different subject classes. For example, teachers could use a theme such as “The Post Office” to get their students to do some writing, to perform some computations, to talk about geography, and even to discuss the difficult position of the working class in other countries. [Karp 2014, p. 315]. 
\\

\leftskip=0cm
His career as a research mathematician was already starting. That same year 1922 he proved his first famous result, in the field of trigonometric series, building an almost everywhere divergent Fourier-Lebesgue series; it was published the following year in the new   \textit{Polishjournal Fundamenta Mathematicæ}\footnote{This journal, founded in 1920 by a group of Polish mathematicians to strenghten the mathematical homeland culture in the years of restoration of Poland independence after the end of the First World War, was at the same time intended with a deep international vocation.}.
\\
\indent
On the  21st January 1924 Lenin died. Under the first years of Joseph Stalin's rule Kolmogorov's career and prestige  took off: he graduated in 1925, and after his postgraduate period, in 1929 he was enrolled in the Moscow University Institute of Mathematics of Mechanics. After ten years, in the dramatic years of the Stalinist terror, Kolmogorov was elected a full member of the USSR Academy of Sciences. In 1933 he had published, in German, his revolutionary treatise on the foundations of probability theory: \textit{Grundbegriffe der Wahrscheinlichkeitsrechnung}, laying the axiomatic foundations of the theory of probability.

\section{Living through the 1930's and 40's: Kolmogorov and science under Stalinism}

Despite the change of regime in the Russian Empire, in the 1920s and early 1930s scientist's in the USSR maintained and developed the ongoing international relationships, specially with Germany, France – a center of emigrés intellectuals after the 1905 upheaval –, and the United States. Legacy of the schools of Dmitry Fëdorovič Egorov\footnote{Dmitry Fëdorovič Egorov (Moscow, 1869- Kazan, 1931) was a Russian mathematician who mainly dealt with differential geometry and mathematical analysis.} and Luzin, Kolmogorov had tight connections in France.\footnote{For more details on relations between France and Russia, see [Deminov 2009].} 
\\
\indent
During a journey with Pavel Alexandrovich Aleksandrov in Germany and France between 1930 and 1931, he met the French mathematicians Maurice Fréchet, Émile Borel and Henri Léon Lebesgue. Deminov asserts that the intense relations between these two countries seemed more flourishing than ever, with no clouds on the horizon: 
\\

\leftskip=1cm
\noindent
Au début des années 30, la recherche mathématique semblait donc pouvoir se poursuivre sans obstacle. La vie normale de la communauté scientifique, interrompue par la première guerre mondiale, avait repris et la collaboration des mathématiciens français avec leurs collègues soviétiques se développait dans les meilleures conditions.\footnote{Eng.tr.: At the beginning of the 1930s, mathematical research seemed to be able to continue without obstacles. The normal life of the scientific community, interrupted by the First World War, had resumed and the collaboration of French mathematicians with their Soviet colleagues was developing under the best conditions.} [Deminov 2009, p 130]
\\

\leftskip=0cm
\noindent
But in the time of the great purges starting at the mid 1930s the connections almost collapsed. As early as 1932, Luzin was denied participation in the International Congress of Mathematicians in Zurich and, two years later, Kolmogorov was unable to go to Paris, despite a scholarship granted by the Rockefeller Foundation. Since the late 1930s, scientist in the Soviet Union hardly traveled abroad. Kolmogorov could not leave the Soviet Union until September 1954, when the Russians were finally able to take part again in the International Congress of Mathematicians in Amsterdam, albeit with a tiny delegation of four people (Alexandrov, Kolmogorov, V. Kozetsky e Sergej Michajlovich Nikolskij).\footnote{Already suspended in 1936, the congress was only reinstated in 1950, but, on that occasion, the entire Russian academic community did not participate in the event. For more details, see in the Proceedings of the ICM 1950, Cambridge, Massachusetts, [Graves, Hille, Smith, Zariski 1955, p 122]}
\\
\indent
During my research work on the cultural origins of the theorem and the research project on classical mechanics and the theory of dynamical systems, which Kolmogorov presented at the Amsterdam congress, I was led to consider his scientific work during the time of terror Stalinist. In fact, he seems to have devoted himself, in a hidden way, to topics of classical mechanics related to Poincaré's research\footnote{see in this regard the preprint [Fascitiello 2022] \textit{Thirty years after: insights on the cultural origins of Andrej N. Kolmogorov's 1954 invariant tori theorem from a short conversation with Vladimir I. Arnold}}, thus continuing an incipient interest in the theory of dynamical systems, in particular in Vito Volterra's seminal studies on population dynamics, as the beginning of a mathematical theoretical biology [Kolmogorov 1935, 1936]. Kolmogorov recalled his early interest in astronomy, thanks to the works of Camille Flamarion (1842–1925), and also in biology:
\\

\leftskip=1cm
\noindent
[...] For a time, interest in other sciences took over. The first big impression on me of the strength and significance of scientific research was made by the book by K. A. Timirjazev\footnote{He refers to Kliment Arkad'evič Timirjazev(1843-1920), a Russian plant physiologist and a major proponent of the Evolution Theory of Charles Darwin in Russia.} The Life of Plants. 
\\

\leftskip=0cm
\noindent
The 23 editions published between 1898 and 1962 in 5 languages - therefore before and after October 1917 - of the famous book by the botanist Timiryazev show his adherence to radical and scientistic ideals and to the vision of science in the Soviet Union. Requirements of Soviet science well described in [Nicolaidis 1990]:
\\

\leftskip=1cm
\noindent
After 1928 however, the Stalinist regime proclaimed a so called "Marxist" official ideological line concerning science. This ideological line became the official line of Soviet astronomy in 1931 Its principles were the following:
\\
\noindent
(1) \textit{There are two sorts of astronomies Soviet and bourgeois}. This principle comes from the dogmatic principle that a capitalist regime restrains the scientific evolution while on the contrary, the construction of the socialist regime implies in addition the construction of a new, superior science. This principle of "two sciences" was the main Stalinist principle concerning all scientific fields. We will see that in astronomy, the application of this principle was to have terrible consequences for the leading Russian astronomers.
\\
\noindent
(2) Soviet astronomy must serve Soviet society more precisely \textit{astronomy must serve ideology and the economy.}
\\
\noindent
But how could astronomy serve Stalinist ideology? [...] Astronomy was a scientific tool that would help to disprove what Stalinists called "religious myths". In a more specifically scientific field, soviet astronomy was ordered to fight against what was termed idealistic western cosmological theories, and especially against the theory of general relativity and the concept of a finite universe - because to put limits and an age to the universe would imply the Creation and so the existence of a God.
\\
\noindent
The relation between astronomy and the Soviet economy was a more complicated concept.
\\
\noindent
The general line that all activities in the USSR must serve the "building of socialism" implied that research in astronomy must also have industrial applications It was difficult to make applications concrete, and so the ideological line spoke about researches concerning Earth Sun relations and geodesy. [Nicolaidis 1990, pp. 346-347]
\\

\leftskip=0cm
\noindent
The risk of researching questions linked to 19th century open problems in celestial mechanics during those years seems very high, taking into account the “purge of astronomers” in years 1936-37 which involved not scientist in various disciplines (celestial mechanics, geology, geophysics, geodesy) working in the wide network of  astronomical observatories in the Russian Empire.\footnote{The causes and evolution of this series of arrestations, emprisonments and murders have been studied in [McCutcheon 1991], [Nicolaidis 1990] and [Eremeeva 1995]; an relevant role was played by the close connection between astronomy and (Orthodox Christian) religion}
\\
\indent
Kolmogorov has often been credited with fully sharing the communist ideals of the Soviet regime. The aim of the present paper is not to fully analize such an issue\footnote{The question of science in the Soviet Union, its impressive devolopment, in the context of an evolution starting in the last fifty years of the Tzarist regime has been paid considerable attention since the dissolution of the USSR and the establishment of the Russian Federation [Kojelnikov 2002], [Gordin et al 2008].}. Nevertheless, there are two opposite interpretations of his behaviour and writings. Loren Graham, in his book \textit{Science in Russia and the Soviet Union. A short History} [Graham 1993], describes Kolmogorov as one of the leading figures in USSR science:
\\

\leftskip=1cm
\noindent
Most people now assume that all influence of Marxism on Soviet science was deleterious. On the contrary, in the works of scientists such as L. S. Vygotsky, A. I. Oparin, V. A. Fock, O. Iu. Schmidt, and A. N. Kolmogorov, the influence of Marxism was subtle and authentic. [Graham 1993, p 3-4].
\\

\leftskip=0cm
\noindent
In his discussion of USSR mathematics, he writes: 
\\

\leftskip=1cm
\noindent
A. N. Kolmogorov, one of Shmidt's authors in the first edition of the \textit{Large Soviet Encyclopedia}, wrote the entry "Mathematics." [...] I will briefly compare Kolmogorov's article with those in the \textit{Encyclopedia Britannica} written by Frank Plumpton Ramsey and Alfred North Whitehead at approximately the same time as the first edition of Kolmogorov's article.
\\
\indent
The points of difference arise on the most essential questions of mathematics: What are the origins of mathematics? and What is the relationship between mathematics and the real world? According to Kolmogorov, mathematics is "the science of quantitative relations and spatial forms of the real world." It arose out of "the most elementary needs of economic life," such as counting objects, surveying land, measuring time, and building structures. In later centuries mathematics became so abstract that its origins in the real world were sometimes forgotten by mathematicians, but Kolmogorov reminded them that "the abstractness of mathematics does not mean its divorce from material reality. In direct connection with the demands of technology and science the fund of knowledge of quantitative relations and spatial forms studied by mathematics constantly grows." Kolmogorov then went on to sketch a history of mathematics in which its growth was intimately related to economic and technological demands. His views were consistent with Lenin's insistence on the material world as the source of human knowledge, and Engels's emphasis on technical needs as a motivating force in the development of knowledge. [Graham 1993, p 118]
\\

\leftskip=0cm
\noindent
In a paper published in 2002 based on his own experiences, George Gunther Lorentz (1910-2006) quotes the testimony and reflection of one of Kolmogorov's most famous students, Vladimir Igorevič Arnold (1937-2010), published in the volume \textit{Kolmogorov in perspective} in 2000 by the American Mathematical Society and the London Mathematical Society [AA.VV. 2000]. Precisely in an attempt to reconstruct the reasons which prompted Kolmogorov in the two-year period 1953-54 to deal with one of the mathematical questions which, after probability, made him more famous - the one which today goes by the name of KAM theory - he writes: 
\\

\leftskip=1cm
\noindent
Although Andrei Nikolaevich himself regarded the hopes that appeared in 1953\footnote{Stalin died in March 1953} as the main stimulus for his work, he always spoke with gratitude about Stalin (following the old principle of saying only nice things about the dead): "First, he gave each academician a quilt in the hard year of the war, and second, he pardoned my fight in the Academy of Sciences, saying, 'such things happen also here'." Andrei Nikolaevich also tried to speak kindly about Lysenko, who had fallen into disfavor, claiming that the latter had sincerely erred out of ignorance (while Lysenko was in power, the relation of Andrei Nikolaevich to this "champion in the struggle against chance in science" was quite different).
\\
\indent
[...] "Some day I will explain everything to you," Andrei Nikolaevich used to tell me after having done something contrary to his principles. Seemingly, pressure was exerted on him by some evil genius whose influence was enormous (the role of the group transmitting the pressure was played by well-known mathematicians). He hardly lived to the times when it became possible to speak of these things, and, like almost all people of his generation who lived through the 1930's and 40's, he was afraid of "them" to his last day. One should not forget that for a professor of that time not to tell the proper authorities about seditious remarks made by an undergraduate or graduate student not infrequently meant being accused the next day of having sympathy with the seditious ideas (in a denouncement by the very same student-provocateur). [Arnold 2000, p. 92]
\\

\leftskip=0cm
\noindent
In the same year that the great purge of Soviet astronomers began, there was an event involving mathematicians and the research group of which Kolmogorov was a part.
\\
\indent
The Luzin affair \footnote{[Levin 1990], [Lorentz 2002],[Katuteladze 2012], [Katuteladze 2013], [Demidov, Lëvshin 2016]}, started with anonymous accusation by the newspaper \textit{Pravda} in eight long articles from July 2 to 16, 1936, as "Enemy under the guise of a Soviet citizen" [Kutateladze 2013, p.A86]. A commission was immediately appointed by the Academy of Sciences and the first session of the process began on July 7, 1936, which was followed by others in a few days. A sentence was reached, which proved to be less drastic than the tones of the articles and the trial, and the Soviet mathematician was acquitted of practically all charges. In this affair was involved Enrst Ko'lman, a Cezch emigrée of Jewish ascendence and Marxist convictions. Aleksey E. Levin writes in [Levin 1990]: 
\\

\leftskip=1cm
\noindent
Kol'man's formal position inside the party hierarchy was head of the mathematical section of the Communist Academy during the first half of the 1930s; he was subsequently promoted early in 1936 to head of the science department of the party's Moscow City Committee (he held this office until 1938 when he apparently lost some powerful protection and was sent to a teaching position). During the 1930s, Kol'man regularly published on the philosophy of mathematics and had many personal connections in mathematical circles, especially among politically active youth.
\\
\noindent
Kol'man's status makes it unlikely that any public attack on Luzin would have been launched without his approval. [Levin 1990, pp 98-99]
\\

\leftskip=0cm
\noindent
A rift opened between Luzin and most of his students, including Kolmogorov and Aleksandrov, who played an active role in the commission set up for the trial.
After the denouement of the matter, the split continued to manifest itself: ten years after the Luzin affair, he voted against the election of Alexandrov to become a full member of the Academy of Sciences:
\\

\leftskip=1cm
\noindent
Lusin died in 1950, but not before a final violent collision with
Aleksandrov and Kolmogorov. In 1946, the Academy had to elect a new group of members, this time with preference to the applied sciences. This allowed Lusin to vote against the topologist Aleksandrov. To everybody’s consternation, as a reaction, Kolmogorov slapped Lusin’s face on the floor of the Academy. The president of the Academy, S. I. Vavilov, was at a loss of what to do. Finally the incident was reported to the Kremlin. It was said that Stalin was not astonished. "This happens even among us," was his reply. In other words, Stalin recommended to do nothing.\footnote{This was the \textit{fight} to which Kolmogorov referred in the words reported by Arnold} [Lorentz 2002, p 207]
\\

\leftskip=0cm

\section{On a new confimation of Mendel’s Laws: Kolmogorov participates in the polemic on genetics in the USSR 1940}

The so-called Lysenko affair is one vexed question on science and ideology in the Soviet Union. A violent campaign against genetics, not considered to conform to dialectical materialism was one of the most devastating political intrusions into Soviet intellectual life under the Stalin regime.\footnote{To name a few, see [Joravsky 1970], [Graham 1993], [Graham 2016], [Ptushenko 2021]} 
\\
\indent
Trochym Denysovych Lysenko (September 29, 1898 Karlivka (Ukraine) - November 20, 1976, Moscow) exerted a growing influence on Russian biology from the mid-1930s onwards, until it reached the vertex of an ascending parabola in the late 1940s, with the approval of Joseph Stalin\footnote{At the Congress of the Communist Party of the Soviet Union in 1935, Lysenko delivered a speech which ended with the applause of the audience and the emblematic exclamation of Joseph Stalin, present at the congress: "Bravo, comrade Lysenko, bravo!"} himself, only to hit rock bottom after more than thirty years, in 1965, when geneticists called him an impostor and attributed to him all the damage caused to Soviet agriculture - just think that in the thirty years 1935-65 Russia, from Europe's granary became an importing country, after a series of failed crops. 
\\
\indent
Among the tenets of what will be termed Lysenkonism are the critique of Mendelism, the denial of the applicability of chemistry, physics and mathematics to the solution of any biological problem [Lysenko 1940], as well as the promise to quickly solve all agricultural tasks set by the party.

\subsection{The Lysenko affair}

I his book \textit{The Lysenko affair} (1970), the American historian of science David Joravsky wrote: 
\\

\leftskip=1cm
\noindent
Thus the Lysenko affair has been pictured as a latter-day version of Galileo versus the Church, or Darwin versus the churches: new science denounced to save old theology. The historical reality was far less high-minded, far more serious. Lysenko's school did not derive from a moribund tradition in science; it rebelled against science altogether. Farming was the basic problem, not theoretical ideology. Not only genetics but all the sciences that impinge on agriculture were tyranically abused by quacks and time-servers for about thirty-five years. The basic motivation was not a dream of human perfectibility but a selfdeceiving arrogance among political bosses, a conviction that they knew better than scientists how to increase farm yields. The Lysenko affair, in short, was thirty-five years of brutal irrationality in the campaign for improved farming, with severe convulsions resulting in the academic disciplines that touch on agriculture. [Joravky 1970, p. vii]
\\

\leftskip=0cm
\noindent
Graham devoted chapter 6 of his book \textit{Science in Russia and the Soviet Union: A Short History} [Graham 1993] to this topic: 
\\

\leftskip=1cm
\noindent
The roots of Lysenkoism lie not in Marxist ideology, but in the social and political context of Soviet Russia in the 1930s. Lysenko originated his ideas outside the circles of Marxist philosophers and outside the community of established geneticists. He was a simple agronomist who developed ideas about plants not very different from those of many practical selectionists of the late nineteenth and early twentieth centuries, but who was able to promote those ideas to an unheralded prominence 	because of the political and social situation in which he found himself. An extremely shrewd but basically uneducated man, he learned how to capitalize on the opportunities that the centralized bureaucracy and ideologically charged intellectual atmosphere presented. Seeing that his ideas would fare better if they were dressed in the garb of dialectical materialism, with the help of a young ideologist he recast his arguments in Marxist terms. 
\\
\noindent
[...] Alarmed that the science of genetics itself might be eclipsed, Vavilov\footnote{It refers to Nikolai Vavilov (November 25, 1887 - January 26, 1943), a Russian agronomist, botanist and geneticist, famous for his expeditions around the world in search of varieties of agricultural plants. He was an exemplary researcher with encyclopedic knowledge: he contributed to genetics, botany, plant physiology, plant breeding, plant systematics and evolution and biochemistry.} abandoned the effort to compromise with Lysenko and pointed out the errors in his biological views. 
\\
\noindent
[...] A few other brave	people continued to speak up against Lysenko. At a conference on genetics in December 1936, A.S. Serebrovskii, an outstanding geneticist and sincere Marxist, called Lysenko's campaign "a fierce attack on the greatest achievements of the twentieth century...an attempt to throw us backward a half-century." [Graham 1993, p. 124, 129-130] 
\\

\leftskip=0cm
\noindent
Lysenko's most famous opponent, Nikolai Vavilov, paid with his life for a heroic attempt to publicly defend the achievements of biology hitherto achieved in the Soviet Union, and to oppose Lysenko's scientifically unsubstantiated ideas. As early as 1935 Lysenko began to attack his colleague, accusing him of having hindered the development of agricultural production in Russia. And, although in 1939 Vavilov was elected president of the VII International Congress of Genetics, it was not enough to maintain his prestige. Accused Vavilov of defending classical Mendelian genetics, considered by party ideologues a "bourgeois pseudoscience", he was imprisoned in 1940 and a year later sentenced to death. He died of starvation in 1943 in the Russian prison of Saratov. He was on of the many biologists arrested, exiled and repressed. 
\\
\indent
In a recent book, Lysenko's ghost (2016) considering recent genetical research, has developed a synthetic but thorough examination of genuine scientifc aspects of the discussion on epigenetics together with the evolution of repression and terror under Stalinism. In the conclusion he writes (pp. 139-140): 
\\

\leftskip=1cm
\noindent
With the realization that the inheritance of acquired characteristics might happen after all, was Lysenko right? No, he was not. Some people may think so because they mistakenly link Lysenko uniquely to the doctrine of acquired characteristics, a belief that has been around for serveral thousand years. Lysenko was a very poor scientist, and the inheritance of acquired characteristics  was actually a small part of what he claimed. 
\\
\indent
The fathers and mothers of epigenetics did not use Lysenko's results but developed their views on the basis of molecular biology. [...] Lysenko disregarded the action of genes [...].
\\
\indent
Does this mean that Lysenko was totally worthless as a practical plant breeder, especially in his early years? No. Lysenko had talents in the field [...]. If Lysenko had lived in a normal democratic country, he would be remembered, if at all, as a talented farmer working away in his fields, employing idiosincratic methods but never garnering much support. None of his methods are employed in Russia today. But in the Soviet Union in the 1930s, a country suffering from famine (caused in large part by the disastrous collectivization effort), the need for quick agricultural remedies was acute, and Lysenko offered them. 
\\

\leftskip=0cm

\subsection{Kolmogorov against Lysenko and Kol'man}

In 1939, the scholar N.I. Ermolaeva published an article entitled (translated into English), \textit{Once again about pea laws}\footnote{In Russian: Yeshche raz o gorokhovikh zakonakh}, which went against the validity of Mendel's principle and where, in particular, it ended that Mendel's principle that self-pollination of hybrid plants resulted in 3: 1\footnote{Mendel's principle, based on the 3:1 ratio, states that the dominant trait is present three times as often as the recessive trait} segregation ratios was false. 
\\
\indent
At that time, Kolmogorov had already long been interested in biology, one of his passions since childhood. The article was brought to Kolmogorov's attention by geneticist Aleksandr Sergeevich Serebrovskii \footnote{[Kolmogorov 1940, p 222]. Serebrovskii (Kusrk, February 18, 1892 - Bolshevo, June 26, 1948) was a prominent Russian geneticist, the founder of the Department of Genetics of Moscow University. His major contributions are due to the genetic study of chicken breeds and the development of poultry farming.}. In a paper published in 1940 \footnote{the episode took place at the turn of the world war - after the sanctioning of the Russo-German pact of August 1939 and a few months before the entry into the war of the USSR after the German invasion of June 1941} in the reports of the USSR Academy of Sciences, he discussed Ermolaeva data as well as T. K. Enin data as discussed by Kol'Man: 
\\

\leftskip=1cm
\noindent
In the discussion on genetics that took place in the autumn of 1939 much attention was paid to checking whether or not Mendel's laws were really true. In the basic	discussion on the validity of the entire concept of Mendel, it was quite reasonable and natural to concentrate on the simplest case, which, according to Mendel, results	in splitting in the ratio 3:1. For this simplest case of crossing $Aa \times Aa$, with the feature A dominating over the feature $a$, it is well known that Mendel's concept leads to the conclusion that in a sufficiently numerous progeny (no matter 	whether it consists of one family or involves many separate families resulting from various pairs of heterozygous parents of type $Aa$) the ratio between the number of 	individuals with the feature A (that is, the individuals of the type $AA$ or $Aa$) to 	the number of individuals with the feature a ($aa$ type) should be dose to the ratio 3:1. T.K. Enin, N.I. Ermolaeva and E. Kol'man have concentrated on checking this 	simplest consequence of Mendel's concept. However, Mendel's concept not only 	results in this simplest conclusion on the approximate ratio $3:1$ but also makes it 	possible to predict the average deviations from this ratio. Owing to this it is the statistical analysis of deviations from the ratio $3:1$ that gives a new, more subtle 	and exhaustive way of proving Mendel's ideas on feature splitting. In this paper we 	will try to indicate what we think to be the most rational methods of such checking 	and to illustrate these methods on the material of the paper by N.I. Ermolaeva. 
\\

\leftskip=0cm
\noindent
And he adds in conclusion a strong attack against Kol'man's contribution: 
\\

\leftskip=1cm
\noindent
Kol'man's paper referred to in the beginning of this note does not contain any new facts; it only analyses Enin's data and is based on a complete misunderstanding of the circumstances set forth in this paper. [Kolmogorov 1940, p. 227]. 
\\

\leftskip=0cm
\noindent
Kolmogorov applied a statistic test, now called the \textit{Kolmogorov} or \textit{Kolmogorov-Smirnov test} in order to analyze in this more polished way the validity of Mendel's ratio. Kolmogorov's paper prompted a reaction by Kol'man and by Lysenko himself. 
\\
\indent
Mathematicians Alan Stark and Eugene Seneta, in the article \textit{A.N. Kolmogorov's defense of Mendelism} published in 2011 [Stark, Seneta 2011], examined Kolmogorov's paper and reused the test used by the Russian mathematician in the data collected by Ermolaeva. They have shown that there were errors not only in Ermolaeva's obviously incorrect calculations - and in the consequent deduction of the inadequacy of Mendel's law - but also in the experiment itself performed by Kolmogorov: 
\\

\leftskip=1cm
\noindent
In the above brief $\chi^2$ analysis we have attempted to use an essentially equivalent test to Kolmogorov’s inasmuch as it relies on the approximate standard normality of the $\Delta$’s, after “cleaning” the data appropriately. So while the conclusion drawn by Kolmogoroff (1940) confirms what is now totally accepted, the evidence in support of this conclusion is not as strong as his paper presents. Of course his statistical technology was well beyond the understanding of Lyssenko (1940) and Kolman (1940), who could hardly argue on the grounds of its incompletely justified application and possible arithmetic error, to data which may have been poorly prepared. Seneta (2004) describes Kolman’s leading role in the attacks on mathematicians and traditional pure mathematics in the Soviet Union during the Stalinist era. 
\\

\leftskip=0cm
\noindent
Futhermore, at the end of the article they write: 
\\

\leftskip=1cm
\noindent
In his defense of Mendelism, Kolmogorov [...] relied simply on data. As we have seen, he ignored the fact that, strictly speaking, his test of him assumed continuous data while the actual data was discrete and in some cases based on inappropriately small numbers [Stark, Seneta 2011, p. 185] 
\\

\leftskip=0cm
\noindent
Kolmogorov had not only corrected the results obtained by Ermolaeva, but his defense against Mendelism first involved Kol'man and Lysenko (even if he was not directly mentioned), as Ermolaeva was a student of him. In fact, Lysenko responded in a comment published in the reports of the Academy:, in \textit{In Response to the Article by A. N. Kolmogorov}, [Lysenko 1940]: 
\\

\leftskip=1cm
\noindent
In "Doklady Akademii Nauk SSSR", Volume XXVII, N 1 of 1940, an article by academician A.N. Kolmogorov "On a new confirmation of Mendel's laws". In this article, the author, wanting to prove the "correctness" and inviolability of Mendel's statistical law, gives a number of mathematical arguments, formulas and even curves. I don't feel competent enough to understand this system of mathematical evidence. Besides, I, as a biologist, don't care whether Mendel was a good or a bad mathematician. I have already published my assessment of Mendel's statistical work several times, stating that he had nothing to do with biology.
In this note I would just like to note that even the above-mentioned article by the famous mathematician A.N. Kolmogorov has nothing to do with biological science. 
\\
\noindent
[...] That is why we biologists do not take the slightest interest in mathematical calculations that confirm the useless statistical formulas of the Mendelists. [Lysenko 1940, p 834-835]
\\

\leftskip=0cm
\indent
Also Kol'man also supported the agronomist's point of view.\footnote{As stated in the note 84 in [Joravsky 1970], p. 414.} 
\\
\indent
It could have ended much worse for the Russian mathematician, but  he was spared from the great purge. There is no doubt that at that time Kolmogorov's fame in the field of probability theory - a field with a strong Russian tradition for more than a century - was undisputed, and this may have favored him over the fates of other academics; moreover the cautious and submissive behavior maintained by Kolmogorov and his friend Aleksandrov\footnote{Kolmogorov and Aleksandrov had lived together since 1935, (ref [Kolmogorov 1986]).} during the period of the Stalinist regime - the same behavior which has led some to believe that they were "friends of the regime" - has to be taken into account. Consider also in [Levin 1990], the description of Kol'man previous attitude to the two mathematical friends: 
\\

\leftskip=1cm
\noindent
When the campaign was over\footnote{He refers to the campaign against Luzin} Kol'man's monograph, \textit{Predmet i metod sovremennoi matematiki}, was published\footnote{Year 1936}. The scholarly qualities of this very primitive and inaccurate book are of no relevance here, although it is worth noting that the author expressed his gratitude to A. N. Kolmogorov and P. S. Aleksandrov for reading a draft manuscript. 
\\

\leftskip=0cm
\noindent
Therefore, only hypotheses, but no certainty about the reasons for his salvation. Nevertheless, quoting Lorentz's words in  [Lorentz 2002, p 183] \textit{The voice of Kolmogorov raised in defense of Mendel’s laws was ignored}.

\section{Final remarks}

What prompted the reaction of Kolmogorov to a paper by a student of Lysenko, who was in those years reaching a great influence? No doubt he was encouraged to do so by a fellow Academician, the genetist Serebrovskii, who could thought that the mathematical authority of Kolmogorov regarding numbers in Mendel's laws of heredity could help their defense of Mendelism and maintain USSR biolgoy in connection with the international accepted ideas. Moreover, the role of mathematical studies in the discussions on evolution was growing in the late 1930s\footnote{See [Kingsland 1985]; [Israel 1993]; [Israel, Millán Gasca 2002]}. At Moscow, the biologist Georgii Frantsevich Gause (1910-1989) has intensely worked on the mathematical analysis of the struggle for life and published two monographies on it in the USA and in France (respectively in 1934 and 1935) . Kolmogorov was sympathetic with this research trend, and this can be linked to his vision of Mathematics as present in the Soviet Encyclopedia. The reaction of Lysenko's against mathematics in biology can be understood as an underlying aspect of the harsh exchange in the Stalinst Russia of 1940.

\section*{Acknowledgement}
I thank Luca Biasco and Ana Millán Gasca for insights and comments on this preprint, which is subject to further improvements.

\newpage

\section*{Bibliography}
\addcontentsline{toc}{section}{Bibliography}
[AA.VV. 2000] \textsc{AA. VV.} (2000) \emph{Kolmogorov in perspective}. Providence, R.I. American Mathematical Society, London Mathematical Society.
\\[1em]
[Abramov 2010] \textsc{Abramov} Alexander, (2010) \emph{Toward a History of Mathematics Education Reform in Soviet Schools
(1960s–1980s)} in [Karp, Vogeli (2010)], pp. 87-140. 
\\[1em]
[Arnold 2000] \textsc{Arnol'd} Vladimir Igorevič (2000) \emph{On A. N. Kolmogorov}. [in AA.VV. \emph{Kolmogorov in perspective}], pp. 89-108.
\\[1em]
[Charpentier, Lesne, Nikolski 2004] \textsc{Charpentier} Eric, \textsc{Lesne} Annick, \textsc{Nikolski}  Nikolai Kapitonovich (eds.) (2004) \emph{L'héritage de Kolmogorov en mathématiques}. Paris, Éditions Belin [eng. tr. Berlin Heidelberg, Springer Verlag 2007].
\\[1em]
[Deminov 2009] \textsc{Deminov} Sergej Sergeevich (2009) \emph{Les relations mathématiques Franco-Russes entre les deux guerres mondiales}. Revue d'histoire des sciences, 2009/1 Tome 62, pp. 119-142.
\\[1em]
[Deminov, Lëvshin 2016] \textsc{Deminov} Sergej Sergeevich, \textsc{Lëvshin} Boris V. (eds.) (2016)\emph{The Case of Academician Nikolai Nikolaevich Luzin}. Providence, Rhode Island, American Mathematical Society.
\\[1em]
[Diner 1992] \textsc{Diner} Simon (1992) \emph{Les voies du chaos déterministe dans l'école russe}. In [Dahan Dalmedico, Chabert, Chemla 1992], pp 331-370.
\\[1em]
[Fascitiello 2022] \textsc{Fascitiello} Isabella (2022) \emph{Thirty years after: insights on the cultural origins of Andrej N. Kolmogorov’s 1954 invariant tori theorem from a short conversation with Vladimir I. Arnold}. Preprint, University of Roma Tre, Department of Mathematics and Physics. URL:
\textit{https://www.researchgate.net/
\\
\noindent
publication/366205118\_Thirty\_years\_after\_insights\_on\_the\_cultural\_origins\_of
\\
\noindent
\_Andrej\_N\_Kolmogorov's\_1954\_invariant\_tori\_theorem\_from\_a\_short\_conversa
\\
\noindent
tion\_with\_Vladimir\_I\_Arnold}
\\[1em]
[Gerretsen, De Groot 1957] \textsc{Gerretsen} Johan C.H., \textsc{De Groot} Johannes (eds.) (1957) \emph{Proceedings of the International Congress of Mathematicians 1954 (Amsterdam September 2 - 9)}. Groningen/Amsterdam, 
Erven P. Noordhoff N.V./North-Holland Publishing Co.
\\[1em]
[Gordin, Hall, Kojevnikov 2008] \textsc{Gordin} Michael, \textsc{Hall} Karl, \textsc{Kojevnikov} Alexei (eds.) (2008) \emph{Intelligentsia Science: The Russian Century, 1860-1960}. Chicago, Chicago University Press, Osiris, 23(1).
\\[1em]
[Graham 1993] \textsc{Graham} Loren R. (1998) \emph{Science in Russia and the Soviet Union: A Short History}. New York, Cambridge University Press.
\\[1em]
[Graham 2016] \textsc{Graham} Loren R. (2016) \emph{Lysenko's ghost. Epigenetics and Russia}. Cambridge/London, Harvard University Press.
\\[1em]
[Graves, Hille, Smith, Zariski 1955] \textsc{Graves} Lawrence M., \textsc{Hille} Einar, \textsc{Smith} Paul A., \textsc{Zariski} Oscar (eds.) (1957) \emph{Proceedings of the International Congress of Mathematicians 1950 (Cambridge, Massachusetts, U.S.A. 1950}. Providence, RI, American Mathematical Society. 
\\[1em]
[Israel 1993] \textsc{Israel} Giorgio (1993) \emph{The emergence of biomathematics and the case of population dynamics: a revival of mechanical reductionism and Darwinism}. Science in context, 6, 469-509.
\\[1em]
[Israel, Millán Gasca 2002], \textsc{Israel} Giorgio, \textsc{Millán Gasca} Ana (2002) \emph{The Biology of numbers. The correspondence of Vito Volterra on mathematical biology}  Basel, Birkhäuser.
\\[1em]
[Joravsky 1970] \textsc{Joravsky} David (1970) \emph{The Lysenko affair}. Chicago/London, University of Chicago Press.
\\[1em]
[Karp 2012] \textsc{Karp} Alexander (2012) \emph{Soviet mathematics education between 1918 and 1931: a time of radical reforms}. Karlsruhe, ZDM Mathematics Education, vol. 44, pp. 551–561.
\\[1em]
[Karp 2014] \textsc{Karp} Alexander, (2014) \emph{Mathematics Education in Russia} in [Karp, Schubring 2014] pp. 303-322.
\\[1em]
[Karp, Schubring 2014] \textsc{Karp} Alexander, \textsc{Schubring} Gert (edited by) (2014) \emph{Handbook on the History of Mathematics Education}. New York, Springer.
\\[1em]
[Karp, Vogeli 2010] \textsc{Karp} Alexander,\textsc{Vogeli} Bruce (eds.) (2010) \emph{Russian mathematics education: history and world significance}. Series of mathematics educations, vol. 4 London/New Jersey/Singapore, World Scientific Publishing Co Pte Ltd.
\\[1em]
[Kingsland 1985] \textsc{Kingsland} Sharon (1985) \emph{Modeling nature: Episodes in the history of population ecology}. Chicago, The Chicago University Press. 
\\[1em]
[Kojevnikov 2002] \textsc{Kojevnikov} Alexei (2002) \emph{Introduction: A New History of Russian Science}. Cambridge University Press. Science in Context 15(2), pp.177–182.
\\[1em]
[Kojevnikov 2008] \textsc{Kojevnikov} Alexei (2008) \emph{The Phenomenon of Soviet Science}. In [Gordin, Hall, Kojevnikov 2008], pp. 115-135.
\\[1em]
[Kolmogorov 1923] \textsc{Kolmogorov} Andrei Nikolaevich (1923) \emph{Une série de Fourier-Lebesgue divergente presque partout}. Warsaw, Fundamenta Mathematicae 4, pp. 324-328.
\\[1em]
[Kolmogorov 1935] \textsc{Kolmogorov} Andrei Nikolaevich (1935) \emph{Deviations from Hardy’s formulas under partial isolation} (in Russian). Doklady Akademii Nauk SSSR, vol 3. pp 129-132
\\[1em]
[Kolmogorov 1936] \textsc{Kolmogorov} Andrei Nikolaevich (1936) \emph{Sulla teoria di Volterra della lotta per l'esistenza}. Giornale dell'Istituto Italiano degli Attuari Attuari, vol 7 pp. 74-80.
\\[1em]
[Kolmogorov 1940] \textsc{Kolmogorov} Andrei Nikolaevich (1940) \emph{On a new confimation of Mendel's Laws} (in Russian). Doklady Akademii Nauk SSSR, 27(1940), 38-42. [English translation in:  \emph{Shiryaev 1992}, pp. 222-227].
\\[1em]
[Kolmogorov 1954] \textsc{Kolmogorov} Andrei Nikolaevich (1991/1957) \emph{The general theory of dynamical systems and classical mechanics}, in [Tikhomirov 1991, vol.1, pp.355-374]. Original edition in Russian in \textsc{Gerretsen Johan C.H., De Groot Johannes}(1957) \emph{Proceedings of the International Congress of Mathematicians 1954 (Amsterdam September 2 - 9)}, North Holland, Amsterdam, vol. 1, pp. 315–333.
\\[1em]
[Kolmogorov 1986] \textsc{Kolmogorov} Andrei Nikolaevich (1986) \emph{Memories of P. S. Aleksandrov}  Russian Mathematical Surveys 41(6), pp. 225-246.
\\[1em]
[Kolmogorov 1988] \textsc{Kolmogorov} Andrei Nikolaevich (1988) \emph{Mathematics: science and profession} (Russian) Moscow, Kvan Library n.64, Nauka. [ita. tr. in \textit{https://sheba.spb.ru/za/kvant64-matprof-1988.htm}]
\\[1em]
[Kutateladze 2012] \textsc{Kutateladze} Semën Samsonovich (2013) \emph{The Tragedy of Mathematics in Russia.}. Siberian Electronic Mathematical Reports, Vol. 9, pp. A85–A100.
\\[1em]
[Kutateladze 2013] \textsc{Kutateladze} Semën Samsonovich (2013) \emph{An epilog to the Luzin case}. Siberian Electronic Mathematical Reports, Vol. 10, pp. A1–A6.
\\[1em]
[Levin 1990] \textsc{Levin}  Aleksey E. (1990) \emph{Anatomy of a Public Campaign: "Academician Luzin's Case" in Soviet Political History}. Slavic Review Vol. 49(1), pp. 90-108
\\[1em]
[Lorentz 1990] \textsc{Lorentz} George Gunter (2002) \emph{Mathematics and Politics in the Soviet Union
from 1928 to 1953}.Journal of Approximation Theory 116, pp. 169–223.
\\[1em]
[Lysenko 1940] \textsc{Lysenko} Trochym Denysovych (1940) \emph{In Response to the Article by A. N. Kolmogorov} ibid., vol.28(9), pp. 832-3.
\\[1em]
[Mazliak 2018] \textsc{Mazliak} Laurent (2018) \emph{The beginnings of the Soviet encyclopedia. The utopia and misery of mathematics in the political turmoil of the 1920s}. Centaurus 60(1-2), pp 25-51.
\\[1em]
[Nicolaidis 1990] \textsc{Nicolaidis} Efthymios (1990) \emph{Astronomy and Politics in Russia in the Early Stalinist Period - 1928-1932}.  Journal for the History of Astronomy 21, pp 345-351.
\\[1em]
[Ptushenko 2021] \textsc{Ptushenko} Vasily V. (2021) \emph{The pushback against state interference in science: how Lysenkoism tried to suppress Genetics and how it was eventually defeated}. Genetics 219(4).
\\[1em]
[Stark, Seneta 2011] \textsc{Stark} Alan, \textsc{Seneta} Eugene (2011) \emph{A.N. Kolmogorov’s defence of Mendelism}. Brazil, Genetics and Molecular Biology vol. 34(2), 177-186.
\\[1em]
[Shiryaev 1989] \textsc{Shiryaev} Albert Nikolayevich (1989) \emph{Kolmogorov: Life and Creative Activities}. The Annals of probability 17(3) pp. 866-944.
\\[1em]
[Shiryaev 1992] \textsc{Shiryaev} Albert Nikolayevich (1992) (ed.) \emph{Selected Works of A. N. Kolmogorov, vol II}. Dordrecht, Springer Science+Business Media.
\\[1em]
[Shiryaev 2000] \textsc{Shiryaev} Albert Nikolayevich (2000) \emph{Andrei Nikolaevich Kolmogorov (April 25, 1903 to October 20, 1987) A Biographical Sketch of His Life and Creative Paths}. in [AA.VV. 2000], pp. 1-88.
\\[1em]
[Tikhomirov 1991] \textsc{Tikhomirov} Vladimir Mikhailovich (1991) \emph{Selected works of A.N. Kolmogorov, vol. I}. Dordrecht, Kluwer Academic Publishers.
\\[1em]
[Timirjazev 1912] \textsc{Timirjazev} Kliment Arkad'evič (1912) \emph{The Life Of The Plant. Ten popular lectures}. London, Foreign Languages Publishing House Moscow. [Authorized English translation who was made by A. Sheremetyeva]
\\[1em]
[Zdravkovska, Duren 2007] \textsc{Zdravkovska} Smilka, \textsc{Duren} Peter Larkin (eds.) (2007) \emph{Golden Years of Moscow Mathematics}.  Second edition, American Mathematical Society, History of Mathematics, vol.6.

\end{document}